\documentclass[12pt]{article}
\usepackage{amsfonts}
\usepackage{amsmath}
\usepackage{mathdots}
\usepackage{graphicx}
\usepackage{hyperref}
\newcommand{\goes}{\rightarrow}
\newcommand{\dis}{\displaystyle}
\newlength{\olcu}
\newlength{\eksi}
\begin{document}
\title{Continuity of multivariate rational functions}
\author{Ali Sinan Sert\"{o}z}
\date{March 2010}
\maketitle

The behavior of a multivariable rational function at its singularities is erratic. The simplest case where we have a chance of understanding its behavior is when the denominator vanishes only at the origin. In two-variable case this rational function defines a surface which either intersects the $z$-axis at one point or wraps around it at the origin. To decide which-happens-when is a tricky process. For this reason not many examples float in the literature.  For example how do we calculate
\[ \lim_{(x,y,z)\goes (0,0,0)}\frac{x^3y^2z}{x^4+y^{12}+z^{14}}, \; \; \text{or} \; \; \lim_{(x,y,z)\goes (0,0,0)}\frac{x^3y^2z^2}{x^4+y^{12}+z^{14}}? \]

For a multivariable rational function whose denominator vanishes only at the origin, the continuity of this function at the origin must certainly be  encoded in the exponents of the variables. The task is therefore to undecipher this code, which is given by the following theorem:

{\bf Theorem: } {\it Let $a_1,\dots,a_N$ be non-negative integers,  $m_1,\dots,m_N$ be positive integers and $c_1,\dots,c_N$ be positive real numbers, where $N>1$. Then
\[ \lim_{(x_1,\dots,x_N)\goes (0,\dots,0)} \;  \frac{x_1^{a_1}\cdots x_N^{a_N}}{c_1x_1^{2m_1}+\cdots +c_Nx_N^{2m_N}} \; \; \text{exists if and only if} \; \; \sum_{i=1}^{N}\frac{a_i}{2m_i}>1. \]
Moreover, when the limit exists, then it is zero. }

{\bf Remarks: } Before we prove this theorem, a few remarks are in order.

$\bullet$ First it is easy to notice that we can take all the $c_i$ as $1$ after re-scaling; define the new coordinates as $X_i=\beta_ix_i$ where $\beta_i>0$ and $\beta_i^{2m_i}=c_i$, $i=1,\dots,N$.
Hence from now on we will take $c_i=1$, $i=1,\dots,N$.

$\bullet$ It is also clear that the only influence of the $a_i$s is to set the rate of growth of the function. Therefore they can be chosen as any non-negative real numbers provided that either we restrict the choice of the variables to  non-negative values or we enter into the realm of complex numbers.

$\bullet$ The $N=1$ case is totally trivial and is slightly different than the general case. In that case the limit exists if and only if $\dis \frac{a_1}{2m_1}\geq 1$. When it exists, the limit is 1 when equality holds and is zero otherwise.

$\bullet$ For notational convenience in the proof, we define $\vec{x}=(x_1,\dots,x_N)$, and set
\[ f(\vec{x})=\frac{\prod_{i=1}^{N}x_i^{a_i}}{\sum_{i=1}^{N}x_i^{2m_i}}. \]
We also define
\begin{eqnarray*}
p &=& \prod_{i=1}^{N} m_i,  \\
p_i &=& p/m_i, \; \; i=1,\dots,N.
\end{eqnarray*}

{\bf Proof of the theorem: } \textbf{First assume that the limit exists.} In this case the limit along any path must also exist and be independent of path. For this purpose set $\dis \lambda=(\lambda_1,\dots,\lambda_N)$ where each $\lambda_i>0$, $i=1,\dots,N$. Restricting $f$ to the path
\[ \vec{x_{\lambda}}(t)=(\lambda_1t^{p_1},\dots,\lambda_Nt^{p_N}), \]
we get
\[ f(\vec{x_{\lambda}}(t))=\left( \frac{\prod_{i=1}^{N} \lambda_i^{a_i}}{\sum_{i=1}^{N} \lambda_i^{2m_i}} \right) \, t^{(a_1p_1+\cdots +a_Np_N)-2p}. \]
As $t\goes 0$, this limit will exist and be independent of $\lambda$ only if the power of $t$ is strictly positive, i.e.
\[ a_1p_1+\cdots a_Np_n-2p>0 \]
or equivalently
\[ \frac{a_1}{2m_1}+\cdots +\frac{a_N}{2m_N} >1,  \tag{*}\label{con} \]
which is precisely the necessary condition we seek.

\textbf{Conversely assume that the inequality (\ref{con}) holds.} We will show that $\dis \lim_{\vec{x}\goes 0}|f(\vec{x})|=0$.

We will use induction on $N$. Clearly there is nothing to prove when $N=1$, since then $f(x_1)=x_1^{a_1-2m_1}$ and (\ref{con}) implies immediately that the required limit exits and is zero.

Now assume $N>1$. Our strategy will be to restrict $|f(\vec{x})|$ to lines parallel to one of the coordinate axes, say the $x_1$-axis, and show that it is bounded along each such line with its maximum value going to zero as the line approaches to the origin.

First we observe that if for some $j$ we have $\dis \frac{a_j}{2m_j}\geq 1$, then
\[ |f(\vec{x})|=|x_1^{a_1}\cdots x_j^{a_j-2m_j} \cdots x_N^{a_n}|\,
\frac{x_j^{2m_j}}{ \sum_{i=1}^{N}x_i^{2m_i}} \leq |x_1^{a_1}\cdots x_j^{a_j-2m_j} \cdots x_N^{a_n}|. \]
By the inequality (\ref{con}), either $a_j-2m_j>0$ or $a_i>0$ for some $i$ other than $j$.
Then  by the sandwich theorem we have $\dis \lim_{\vec{x}\goes 0}|f(\vec{x})|=0$.

Therefore we are reduced to the case where $0\leq a_i<2m_i$, $i=1,\dots,N$. It is clear that when (\ref{con}) holds, at least one of the $a_i$ is strictly positive. Without loss of generality assume that $0<a_1<2m_1$.

At this point we quote our induction hypothesis:
\[ \text{If} \;  \frac{d_2}{2m_2}+\cdots+ \frac{d_N}{2m_N}>1,  \text{then} \; \lim_{(x_2,\dots,x_N)\goes(0,\dots,0)} \; \frac{\prod_{i=2}^{N}|x_i|^{d_i}}{\sum_{i=2}^{N}x_i^{2m_i}}=0, \]
where $d_2,\dots,d_N$ are non-negative integers, and $m_2,\dots,m_N$ are positive integers.

Now for any $\vec{x}=(x_1,\dots,x_N)$ set $\pi(\vec{x})=(|x_2|,\dots,|x_N|)$.

We fix $\vec{x}$ and consider the non-trivial case when $\pi(\vec{x})\not = (0,\dots,0)$.

We now restrict the function $f(\vec{x})$ to the line
\[ t\mapsto (t,|x_2|,\dots,|x_N|), \; t\in [0,\infty). \]
Call the restriction of $f$ to this line by $\phi_{\pi(\vec{x})}$;
\[ \phi_{\pi(\vec{x})}(t)=f(t,|x_2|,\dots,|x_N|)=\left( \prod_{i=2}^{N}|x_i|^{a_i}\right) \,
\frac{t^{a_1}}{t^{2m_1}+\left( \sum_{i=2}^{N}x_i^{2m_i}\right)}, \; t\in [0,\infty). \]

Clearly $\dis \phi_{\pi(\vec{x})}(t)\geq 0$ on its domain, $\dis \phi_{\pi(\vec{x})}(0)=0$ and moreover $\dis \lim_{t\goes\infty}\phi_{\pi(\vec{x})}(t)=0$. 
Hence the function $\phi_{\pi(\vec{x})}(t)$ will attain its maximum value at some point, say $t_{\pi(\vec{x})}\in [0,\infty)$. We then have
\[ 0\leq |f(\vec{x})|=\phi_{\pi(\vec{x})}(|x_1|)\leq \phi_{\pi(\vec{x})}(t_{\pi(\vec{x})}), \text{  for all } |x_1|\in [0,\infty).  \]
It now remains to show that $\dis \lim_{\pi(\vec{x})\goes 0}\phi_{\pi(\vec{x})}(t_{\pi(\vec{x})})=0$.

A direct calculation yields that $\phi_{\pi(\vec{x})}(t)$ has its maximum at
\[  t_{\pi(\vec{x})}=\left( \frac{a_1}{2m_1-a_1}\right)^{\frac{1}{2m_1}} \left( \sum_{i=2}^{N} x_i^{2m_i}\right)^{\frac{1}{2m_1}}. \]
The maximum value of $\phi_{\pi(\vec{x})}(t)$ can now be written as
\[ \phi_{\pi(\vec{x})}(t_{\pi(\vec{x})})= K \, g(\pi(\vec{x}))^{(1-\frac{a_1}{2m_1})}, \]
where $K$ is a constant and
\[ g(\pi(\vec{x}))=\frac{\prod_{i=2}^{N}|x_i|^{d_i}}{\sum_{i=2}^{N}x_i^{2m_i}}, \]
where $\dis d_i=\frac{a_i}{1-\frac{a_1}{2m_1}}$, $i=2,\dots,N$. (Compare this with our induction hypothesis above.)

The condition (\ref{con}) implies that
\[ \frac{d_2}{2m_2}+\cdots+ \frac{d_N}{2m_N}=\left( \frac{1}{1-\frac{a_1}{2m_1}}\right) \left( \frac{a_2}{2m_2}+\cdots+ \frac{a_N}{2m_N}\right) > 1 \]
and this in turn, by the induction hypothesis, implies that
\[ \lim_{\pi(\vec{x})\goes 0}\phi_{\pi(\vec{x})}(t_{\pi(\vec{x})})=0, \]
which completes the proof. \hfill $\Box$

We can  discuss even the differentiability of such fractions:\\
{\bf Corollary: } {\it Let $a_1,\dots,a_N,m_1,\dots,m_N$ be all positive integers and $c_1,\dots,c_N$ be positive real numbers, where $N>1$. Then the function
\[ f(\vec{x})=\frac{\prod_{i=1}^{N}x_i^{a_i}}{\sum_{i=1}^{N}c_ix_i^{2m_i}} \]
is $C^1$ at the origin if
\[ \sum_{i=1}^{N}\frac{a_i}{2m_i}>1+\max_{1\leq j\leq N}\{ \frac{a_j}{2m_j}\} . \] }

{\bf Proof: } We calculate the $j$th partial derivative for $j=1,\dots,N$ and find that
\[ \left| \frac{\partial f}{\partial x_j} \right| \leq \frac{|x_j^{a_j-1}|\,
\prod_{ i=1 , i\not = j}^{N}|x_i|^{a_i}}{\sum_{i=1}^{N}c_ix_i^{2m_i}} \left( |a_j-2m_j|+|a_j|\right). \]
Now apply the theorem to assure the continuity of this expression at the origin.
\hfill $\Box$

\textbf{A final remark:} The proof of the theorem reveals that there is a distinguished path, $(\lambda_1t^{p_1},\dots,\lambda_Nt^{p_N})$, with the property that the  limit exits if and only if it exists along this  path. It is tempting to ask at this point if such a \emph{royal path} exist for every limit problem.

{\it Bilkent University, Department of Mathematics, 06800 Ankara, Turkey.} \\
{\tt sertoz@bilkent.edu.tr} \\
See also this  \href{http://sertoz.bilkent.edu.tr/depo/sertoz_theorem.htm}{link} for some humor.

\end{document}